\documentclass{siamltex}
\usepackage{amssymb,amsmath,color,bm,bbm,mathrsfs,amscd}
\usepackage{dsfont}
\usepackage{tikz}
\graphicspath{{fig/}}

\def\N{\mathbb{N}}
\def\R{\mathbb{R}}

\DeclareMathOperator*{\argmin}{argmin}
\DeclareMathOperator*{\sgn}{sgn}

\newcommand{\real}{{\mathbb{R}}}
\newcommand{\realnoneg}{\R_+}

\newcommand{\ov}{\overline}

\newcommand{\map}[3]{#1: #2 \rightarrow #3}

\newcommand{\se}{\text{ if }}
\newcommand{\summ}{\sum\limits}

\newcommand{\1}{\mathbbm{1}}

\newcommand{\mc}{\mathcal}

\newcommand{\be}{\begin{equation}}
\newcommand{\ee}{\end{equation}}
\newcommand{\ba}{\begin{array}}
\newcommand{\ea}{\end{array}}
\renewcommand{\l}{\left}
\renewcommand{\r}{\right}
\newcommand{\ds}{\displaystyle}
\newcommand{\de}{\mathrm{d}}

\newcommand{\simplex}{\Pi}
\newcommand{\graph}{\mc G}
\newcommand{\nodeset}{\mc V}
\newcommand{\linkset}{\mc E}

\newcommand{\ksmargin}[1]{\marginpar{\color{red}\footnotesize [KS]:#1}}
\renewcommand{\ksmargin}[1]{}

\newcommand{\delayfunc}{T}

\newcommand{\paths}{\mc P}
\newcommand{\path}{p}

\newcommand{\interior}{\mathrm{int}}

\newcommand{\onebf}{\mathbf{1}}

\newcommand{\incidencemat}{A}

\newcommand{\edgeset}{\mc E}

\newcommand{\boundconst}{\kappa}
\def\qed{\hfill \vrule height 7pt width 7pt depth 0pt              \medskip}
\newtheorem{remark}{Remark}
\newtheorem{assumption}{Assumption}
{\newtheorem{example}{Example}}

\title{Stability analysis of transportation networks with multiscale driver decisions
\thanks{
This work was supported in part by NSF EFRI-ARES
grant number 0735956. Any opinions, findings, and conclusions
or recommendations expressed in this publication are
those of the authors and do not necessarily reflect the views
of the supporting organizations.}}
\author{Giacomo Como
\thanks{Laboratory for Information and Decision Systems,
Massachusetts Institute of Technology, 77 Mass Ave, Cambridge (MA), 02139, US
({\tt giacomo@mit.edu})}.
\and Ketan Savla
\thanks{Laboratory for Information and Decision Systems,
Massachusetts Institute of Technology, 77 Mass Ave, Cambridge (MA), 02139, US
({\tt ksavla@mit.edu})}.
\and Daron Acemoglu
\thanks{Department of Economics,
Massachusetts Institute of Technology, 77 Mass Ave, Cambridge (MA), 02139, US
({\tt daron@mit.edu})}.
\and Munther A.~Dahleh
\thanks{Laboratory for Information and Decision Systems,
Massachusetts Institute of Technology, 77 Mass Ave, Cambridge (MA), 02139, US
({\tt dahleh@mit.edu})}. 
\and Emilio Frazzoli 
\thanks{Laboratory for Information and Decision Systems,
Massachusetts Institute of Technology, 77 Mass Ave, Cambridge (MA), 02139, US
({\tt frazzoli@mit.edu})}. 
} \date{}

\begin{document}

\maketitle

\begin{abstract}
Stability of Wardrop equilibria is analyzed for dynamical transportation networks in which the drivers' route choices are influenced by information at multiple temporal and spatial scales. The considered model involves a continuum of indistinguishable drivers commuting between a common origin/destination pair in an acyclic transportation network. The drivers' route choices are affected by their, relatively infrequent, perturbed best responses to global information about the current network congestion levels, as well as their instantaneous local observation of the immediate surroundings as they transit through the network. A novel model is proposed for the drivers' route choice behavior, exhibiting local consistency with their preference toward globally less congested paths as well as myopic decisions in favor of locally less congested paths. The simultaneous evolution of the traffic congestion on the network and of the aggregate path preference is modeled by a system of coupled ordinary differential equations.  The main result shows that, if the frequency of updates of path preferences is sufficiently small as compared to the frequency of the traffic flow dynamics, then the state of the transportation network ultimately approaches a neighborhood of the Wardrop equilibrium. The presented results may be read as a further evidence in support of Wardrop's postulate of equilibrium, showing robustness of it with respect to non-persistent perturbations. The proposed analysis combines techniques from singular perturbation theory, evolutionary game theory, and cooperative dynamical systems. 

\end{abstract}

\section{Introduction}
As transportation demand is dramatically approaching its infrastructure capacity, a rigorous understanding of the relationship between the macroscopic properties of transportation networks and realistic driver route choice behavior is attracting renewed research interest. Such an analysis is essential, among other things, for appropriate design of incentives influencing drivers' behavior in order to induce a desired socially optimal usage of the transportation infrastructure. 
A particularly relevant issue is the impact of drivers' {\it en route} responses to \emph{unexpected events} on the overall transportation network dynamics. This issue is particularly significant in a modern real-life transportation network scenario, where recent technological advancements in intelligent traveller information devices have enabled drivers to be much more flexible in selecting their routes to destination even while being {\it en route}. While there has been a significant research effort to investigate the effect of such technologies on the route choice behavior of drivers, e.g., see \cite{Srinivasan.Mahmassani:00,Polydoropoulou.BenAkiva.ea:96}, the analytical study of the dynamical properties of the whole network under such behavior has attracted very little attention.

This paper is focused on the stability analysis of transportation networks in a setup where the drivers have access to traffic information at multiple temporal and spatial scales and they have the flexibility to switch their route to destination at every intermediate traffic intersection. Specifically, we consider a model in which the drivers choose their routes while having access to relatively infrequent global information about the network congestion state, and real-time local information as they transit through the network. The drivers' route choice behavior is then influenced by relatively slowly evolving path preferences as well as myopic responses to the instantaneous observation of the local congestion levels at the intersections. This setup captures many real-life scenarios where unexpected events observed {\it en route} might cause drivers to take a temporary detour, but not necessarily to change their path preferences. Such path preferences may instead be updated, e.g., on a daily, weekly, or longer time basis, in response to information about the global congestion state of the different origin-destination paths collected from the drivers' personal experience, their opinion exchanges with their peers, as well as from information media. However, since the traffic dynamics is significantly influenced by the drivers' response to real-time local information, such responses can influence the drivers' path preference thereby modifying their global route choice behavior in the long run. 
We propose and analyze a novel model for the drivers' route choice behavior that combines relatively infrequent information about the global congestion status of the network with real-time local observations as explained below. 

In our model, the network is represented by a directed acyclic graph with one origin and one destination. A continuous constant flow of indistinguishable drivers enters from the origin, and flows through the network until reaching the destination node. Traffic parameters, such as average speed, traffic density, and flow, are modeled as homogeneous quantities on every link, related one to each other by functional dependencies representative of the links' congestion properties. The dynamics of such traffic parameters is governed by the law of conservation of mass, as well as the drivers' route choice behavior. In turn, the drivers' route choice behavior is assumed to be influenced by two factors: the aggregate path preference, measuring the relative appeal of the different routes to the drivers, and local observations of the current congestion levels. The path preference dynamics evolve at a slow time scale (as compared to the traffic dynamics), following a perturbed best response to global information, embodied by the current congestion levels on the whole network. When traversing an intermediate node in the network, drivers behave according to their path preference, if this is consistent with the current, locally observed, aggregate behavior of the other drivers. On the other hand, when there is a discrepancy between the aggregate path preference and the locally observed aggregate behavior, then drivers tend to compensate this by myopically preferring routes which appear to be locally less congested. 

The above-described model gives rise to a double feedback dynamics, governed by a finite-dimensional system of coupled ordinary differential equations. Such a dynamical system has two natural time scales, characterizing the dynamics of the drivers' aggregate path preference and of the traffic parameters on the different links, respectively. We study the long-time behavior of this dynamical system: our main result shows that, in the limit of small update rate of the aggregate path preferences, a state of approximate \emph{Wardrop equilibrium} \cite{Wardrop:52} is approached. The latter is a configuration in which the delay associated to any source-destination path chosen by a nonzero fraction of the drivers does not exceed the delay associated to any other path. Our results contribute to providing a stronger evidence in support of the significance of Wardrop's postulate of equilibrium for a transportation network. They may also be read as a sort of robustness of such equilibrium notion with respect to non-persistent perturbations of the network.   

The analytical arguments we propose mainly rely on three ideas: adopting a singular perturbation approach \cite{Khalil:96}, by considering the aggregate path preference as `quasi-static' when studying the fast scale dynamics of the traffic parameters, and the traffic parameters as `almost equilibrated' when analyzing the slow scale dynamics of the aggregate path preference; exploiting the inherent cooperative\footnote{Here, the adjective `cooperative' is intended in the sense of Hirsch \cite{Hirsch:82,Hirsch:85}.} dependence of the route choice function on the local traffic parameters in order to establish exponential stability of the fast scale dynamics of the traffic parameters; adapting results from evolutionary population games \cite{HofbauerSigmund:03,Sandholm:10} in order to establish stability properties of the slow scale perturbed best response dynamics of the aggregate path preference.  
 
%
Our work is naturally related to two streams of literature on transportation networks. On the one hand, traffic flows on networks have been widely analyzed with fluid-dynamical and kinetic models: see, e.g., \cite{Garavello.Piccoli:06}, and references therein. As compared to these models (typically described by  partial, or integro-differential equations), ours significantly simplifies the evolution of the traffic parameters (treating them as homogeneous quantities on the links, representative of spatial averages), whereas it highlights the role of the drivers' route choice behavior with its double feedback dynamics, which is typically neglected in that literature. 

On the other hand, transportation networks have been studied from a decision-theoretic perspective within the framework of \emph{congestion games}~\cite{Beckmann.McGuire.ea:56,Rosenthal:73}. In these models, drivers make sequential myopic route choice decisions in pursuit of minimizing their personal travel times, in response to complete information about the whole network. Congestion games are known to belong to a class of games known as \emph{potential games}, a consequence of which is that, best responses of the drivers are \emph{aligned} with the gradient of a common potential function and hence the system eventually converges to a critical point of this potential function, which, under appropriate monotonicity conditions of the congestion properties of the links of the network, corresponds to a Wardrop equilibrium. Such an approach has been used, for example in \cite{Marden.Arslan.ea:TAC09}. The stability of Wardrop equilibrium in the context of communication networks has been studied in \cite{Borkar.Kumar:03}. It is important to note that the two salient features of a typical congestion game setup are that information is available to the drivers at a single temporal and spatial scale, and that the dynamics of traffic parameters are completely neglected by assuming that they are instantaneously equilibrated. In contrast, we study the stability of Wardrop equilibrium in a setting where the dynamics of the traffic parameters are not neglected, and the drivers' route choice decisions are affected by, relatively infrequent global information, as well as their real-time local information as they transit through the network. As a consequence, classic results of evolutionary game theory and population dynamics \cite{HofbauerSigmund:03,Sandholm:10} are not directly applicable to our framework, and novel analytical tools have to be developed, particularly for the analysis of the fast scale dynamics of the traffic parameters.

The rest of the paper is organized as follows. 
In Section~\ref{sec:formulation}, we formulate the model and state the main result. 
%
Section~\ref{sec:proofs} is a technical section that contains the proofs for the main result including intermediate results.
In Section~\ref{sec:sim}, we report results from illustrative numerical experiments.
Finally, we conclude in Section~\ref{sec:conclusion} and also mention potential future research directions.

Before proceeding, we establish here some notation to be used throughout the paper. 
Let $\real$ be the set of reals, $\realnoneg:=\{x\in\R:\,x\ge0\}$ be the set of nonnegative reals.
Let $\mc A$ and $\mc B$ be finite sets. Then, $|\mc A|$ will denote the cardinality of $\mc A$, $\R^{\mc A}$ (respectively, $\R_+^{\mc A}$) the  space of real-valued (nonnegative-real-valued) vectors whose components are labeled by elements of $\mc A$, and $\R^{\mc A\times\mc B}$ the space of matrices whose real entries labeled by pairs of elements in $\mc A\times\mc B$. 
The transpose of a matrix $M \in \real^{\mc A \times\mc B}$, will be denoted by $M'\in\R^{\mc B\times\mc A}$, while $I$ to be an identity matrix, and $\onebf$ the al-one vector, whose size will be clear from the context.
The simplex of probability vectors over $\mc A$ will be denoted by 
$\mc S(\mc A):=\{x\in\R_+^{\mc A}:\,\sum_{a\in\mc A}x_a=1\}$.
If $\mc B\subseteq\mc A$, $\1_{\mc B}:\mc A\to\{0,1\}$ will stand for the indicator function of $\mc B$, with $\1_{\mc B}(a)=1$ if $a \in B$, $\1_{\mc B}(a)=0$ if $a\in\mc A\setminus\mc B$. 
For $p \in [1,\infty]$, $\|\,\cdot\,\|_p$ is the $p$-norm. By default, let $\|\cdot\|:=\|\,\cdot\,\|_2$ denote the Euclidean norm.  Let $\interior(\mc X)$ be the interior of a set $\mc X\subseteq\R^d$, and $\partial \mc X$ denote its boundary. Let $\sgn:\R\to\{-1,0,1\}$ be the sign function, defined by $\sgn(x)$ is 1 if $x>0$, $\sgn(x)=-1$ if $x<1$, and $\sgn(x)=0$ if $x=0$. Conventionally, we shall assume the identity $\de|x|/\de x=\sgn(x)$ to be valid for every $x\in\R$, including $x=0$. 


\section{Model formulation and main result}
\label{sec:formulation}
In this section, we formulate the problem and state the main result. 
In our formulation, we represent the dynamics of the traffic and the route choice behavior on a transportation network as a system of coupled ordinary differential equations with two time scales representative of route choice behavior influenced by the two levels of information. The key components of our model are: network topology, congestion properties of the links, path preference dynamics, and node-wise route choice decision. We next describe these components in detail.


\subsection{Network characteristics}\label{sect:networkcharacteristics}
Let the topology of the transportation network be described by a directed graph (shortly, \emph{di-graph}) $\graph=(\nodeset,\linkset)$, where $\nodeset$ is a finite set of nodes and $\linkset\subseteq\mc V\times\mc V$ is the set of (directed) links. For every node $v\in\mc V$, we shall denote by $\mc E^-_v$, and $\mc E^+_v$, the sets of its incoming, and, respectively, outgoing links. 
A length-$l$ (directed) path from $u\in\mc V$ to $v\in\mc V$ is an $l$-tuple of consecutive links $\{(v_{j-1},v_j)\in\mc E:\,1\le j\le l\}$ with $v_0=u$, and $v_l=v$. A cycle is path of length $l\ge1$ from a node $v$ to itself. 
Throughout this paper, we shall assume that:  
\begin{assumption}\label{assumption1}
The di-graph $\mc G$ contains no cycles, has a unique origin (i.e., some $v\in\mc V$ such that $\mc E^-_v=\emptyset$), and a unique destination (i.e., $v\in\mc V$ such that $\mc E^+_v=\emptyset$). Moreover, there exists a path to the destination node from every other node in $\mc V$. 
\end{assumption}
\\ Assumption \ref{assumption1} implies that one can find a (not necessarily unique) topological ordering of the node set $\mc V$ (see, e.g., \cite{Cormen.Leiserson:90}). We shall assume to have fixed one such ordering, identifying $\mc V$ with the integer set $\{0,1,\ldots,n\}$, where $n:=|\mc V|-1$, in such a way that 
$$\mc E^-_{v}\subseteq\bigcup\nolimits_{0\le u<v}\mc E^+_{u}\,,\qquad\forall  v=0,\ldots,n\,.$$

We shall model the traffic parameters as time-varying quantities which are homogeneous over each link of the network. 
Specifically, for every link $e\in\mc E$, and time instant $t\ge0$, we shall denote the current traffic density, and flow, by $\rho_e(t)$, and $f_e(t)$, respectively, while $$\rho(t):=\{\rho_e(t):\,e\in\mc E\}\,,\qquad f(t):=\{f_e(t):e\in\mc E\}$$ 
will stand for the vectors of all traffic densities, and flows, respectively. 
Current traffic flow and density on each link are related by a functional dependence 
\be\label{flowdensityeqn}f_e=\mu_e(\rho_e)\,,\qquad e\in\mc E\,.\ee 
Such functional dependence models the drivers' speed and lane adjustment behavior in response to traffic density on a particular segment of a road.
It will be assumed to satisfy the following: 
\begin{assumption}\label{assumption2}
For every link $e \in \edgeset$, the flow-density function $\mu_e:\R_+\to\R_+$ is continuously differentiable, strictly increasing, strictly concave 
and is such that $\mu_e(0)=0$ and $\lim_{\rho_e\downarrow 0}\de \mu_e/\de \rho_e(\rho_e)<+\infty$. 
\end{assumption}
\begin{remark}\label{remarkflow}
Flow-density functions commonly used in transportation theory typically are not globally increasing, but rather have a $\cap$-shaped graph \cite{Garavello.Piccoli:06}: $\mu_e(\rho_e)$ increases from $\mu_e(0)=0$ until achieving a maximum $C_e=\mu_e(\tilde \rho_e)$, and then decreases for $\rho_e\ge \tilde \rho_e$. Assumption \ref{assumption2} remains a good approximation of this setting, provided that $\rho_e$ stays in the interval $[0,\tilde \rho_e)$.
\end{remark}

\noindent For every link $e\in\mc E$, let 
$$C_e:=\sup\{\mu_e(\rho_e):\,\rho_e\ge0\}=\lim_{\rho_e\to+\infty}\mu_e(\rho_e)$$
be its \emph{maximum flow capacity}. Moreover, let 
$$\mc F_v:=\times_{e\in\mc E^+_v}[0,C_e)\,,\qquad\mc F:=\times_{e\in\mc E}[0,C_e)$$ 
be the sets of local, and, respectively, global admissible flow vectors. 
Observe that our formulation allows for both the cases of bounded and unbounded maximum flow capacities. As the flow $f_e$ is the product of speed and density, it is natural to introduce the delay function
\begin{equation}
\label{eq:delay-function}
\map{T}{\R_+^{\mc E}}{[0,+\infty]^{\mc E}}\,,\qquad
\delayfunc_e(f_e):= \left\{\begin{array}{lcl} 
    +\infty&\se& f_e\ge C_e\\[5pt]
    \mu_e^{-1}(f_e)/f_e& \se & f_e\in(0,C_e) ,\\[5pt]
   1/\frac{\de \mu_e}{\de \rho_e}(0) &\se& f_e = 0\,,
\end{array}  \right. \quad\end{equation}
whose components measure the flow-dependent time taken to traverse the different links.\footnote{Here it has implicitly been assumed, without any loss of generality, that all the links are of unit length.} 
 \begin{example} 
A flow-density function that satisfies Assumption \ref{assumption2} is given by
\begin{equation}
\label{eq:flowfunc-example}
\mu_e(\rho_e) =  C_e\left(1 - e^{- \theta_e \rho_e}\right) \quad \forall e \in \edgeset,
\end{equation}
where $C_e>0$, and $\theta_e>0$. The corresponding delay function is 
$$T_e(f_e)=\frac1{\theta_ef_e}\log\frac{C_e}{C_e-f_e}\,.$$
\end{example}

We shall denote by $\mc P$ the set of distinct paths in $\graph$ from the origin node $0$ to the destination node $n$, and let
$$\incidencemat \in \R^{\mc E \times \mc P}\,,\qquad\incidencemat_{ep}=\left\{\ba{lcl}1&\text{ if }&e\in p\\0&\text{ if }&e\notin p\,,\ea\r.$$
be the link-path incidence matrix of $\graph$. 
The relative appeal of the different paths to the drivers will be modeled by a time-varying probability vector over $\mc P$, which will be referred to as the current \emph{aggregate path preference}, and denoted by $\pi(t)$. 
If one assumes, as we shall do throughout this paper, a constant unit incoming flow in the origin node, it is natural to consider the vector 
$$f^{\pi}:=A\pi$$ 
of the flows associated to the current aggregate path preference. 
Indeed, $f^{\pi}_e=\sum_pA_{ep}\pi_p$ represents the total traffic flow that a link $e\in\mc E$ would sustain in a hypothetic equilibrium condition in which the fraction of drivers choosing any path $p\in\mc P$ is given by $\pi_p$.
Now, let 
$$\Pi:=\left\{\pi\in\mc S(\mc P)\,:\,(A\pi)_e< C_e,\,\forall e\in\mc E\right\}$$
be the set of \emph{feasible path preferences}. Here, the term `feasible' refers to the fact that the flow vector $f^{\pi}$ associated to any $\pi\in\Pi$ satisfies the capacity constraint $f^{\pi}_e<C_e$ for every $e\in\mc E$.  
Observe that, whenever $C_e>1$ for every $e\in\mc E$ (or, in particular, when link capacities are infinite), the set of admissible path preferences $\Pi$ coincides with the whole simplex $\mc S(\mc P)$. 
In contrast, when $C_e\le 1$ for some $e\in\mc E$, $\Pi\subset\mc S(\mc P)$ is a strict inclusion. 
On the other hand, whether $\Pi$ is empty or not depends solely on the value of the min-cut capacity of the network \cite[Ch.~4]{Ahuja.Magnanti.ea:93}
$$C^{*}:=\min_{\substack{\mc U\subseteq\mc V:\\0\in\mc U,n\notin\mc U}}C_{\mc U}\,,\qquad C_{\mc U}:=\summ_{\substack{e=(u,v)\in\mc E\,:\\u\in\mc U,\,v\in\mc V\setminus\mc U}}C_{e}\,,$$
as shown in the following, simply established, result. 
\begin{proposition}
The set $\Pi$ is nonempty if and only if $C^{*}>1$.
\end{proposition}
\proof
Fix a cut-set $\mc U\subseteq\mc V$ such that $0\in\mc U$, and $n\notin\mc U$. Then, every path $p\in\mc P$ contains exactly one link $(u,v)\in p$ such that $u\in\mc U$, and $v\in\mc V\setminus\mc U$. Hence, for every $\pi\in \Pi$, one has that 
$$C_{\mc U}=\summ_{\substack{e=(u,v)\in\mc E:\\u\in\mc U,v\in\mc V\setminus\mc U}}C_e
>\summ_{p\in\mc P}\summ_{\substack{e=(u,v)\in\mc E:\\u\in\mc U,v\in\mc V\setminus\mc U}}A_{ep}\pi_p
=\summ_{p\in\mc P}\pi_p=1\,.$$
Minimizing over all cut-sets $\mc U$ shows that $C^{*}>1$ is necessary for $\Pi$ to be nonempty. 

For the inverse implication, consider a network with the same topology $\mc G$ and link capacities $c_e=\max\{C_e-|\mc E|^{-1}(C^{*}-1),0\}$. 
The min-cut capacity of this network satisfies $c^{*}\ge C^{*}-(C^{*}-1)=1$. Note also that, from our construction, $C_e > c_e \geq 0$.
Therefore, the max-flow min-cut theorem (see, e.g., \cite[Thm.~4.1]{Ahuja.Magnanti.ea:93}) implies that there exists some $\pi\in\Pi$, thus proving that $\Pi$ is nonempty. 
\qed

In the case when $C^{*}\le 1$ it is not hard to show that, since the incoming flow exceeds the outgoing flow of the network, the system will grow unstable, i.e., $\rho_e(t)$ is unbounded as $t$ grows large, for some link $e\in\mc E$. Therefore, throughout this paper we shall confine ourselves to transportation networks satisfying:
\begin{assumption}\label{assumption2.5}
The min-cut capacity satisfies
$C^{*}>1$.
\end{assumption}

\subsection{Route choice behavior and traffic dynamics}
We now describe the drivers' route choice behavior and traffic dynamics on the network.  
We envision a continuum of indistinguishable drivers traveling through the network. 
Drivers enter the network from the origin node $0$ at a constant unit rate, travel through it, and leave the network from the destination node $n$. 
While inside the network, drivers occupy some link $e\in\mc E$. 
The time required by the drivers to traverse link $e$, and the current flow on such link are governed by its congestion properties, as given by \eqref{eq:delay-function}, and (\ref{flowdensityeqn}), respectively. 
When entering the network from the origin node $v=0$, as well as when reaching the tail node $v\in\{1,2,\ldots,n-1\}$ of some link $e\notin\mc E^-_n$, 
the drivers instantaneously join some link $e\in\mc E^+_v$. 
In this paper, we shall model the choice of such new link to depend on infrequently updated perturbed best responses of the drivers 
to global information about the congestion status of the whole network as well as on their instantaneous observation of the local congestion levels. 
We next describe these two aspects of the model in detail.


\subsubsection*{Aggregate path preference dynamics}
The drivers' aggregate path preference $\pi(t)$, already introduced in Sect.~\ref{sect:networkcharacteristics}, models the relative appeal of the different paths to the drivers' population. 
%
The aggregate path preference $\pi(t)$ is updated as drivers access global information about the current congestion status of the whole network.
This occurs at some rate $\eta>0$, which will be assumed small with respect to the time-scale of the network flow dynamics. 
Information about the current status of the network is embodied by the current traffic flow vector $f(t)$. 
From $f(t)$, drivers can evaluate the vector $A'T(f(t))$, whose $p$-th component $\sum_{e \in \mc E} A_{ep}T_e(f_e(t))$ coincides with the total delay one expects to incur on path $p$ assuming that the congestion levels on such path won't change. 
Drivers' are assumed to react to such global information by a \emph{perturbed best response} 
\be\label{Fhdef}F^h(f):=\argmin_{\omega \in\simplex_h } \big\{ \omega'\incidencemat'T(f)  +h(\omega) \big \}\,,\ee 
where $h:\Pi_h\to\R$ is an \emph{admissible perturbation}, satisfying the following: 
\begin{assumption} \label{assumption3}
An admissible perturbation is a function $\map{h}{\Pi_h}{\R}$ where $\Pi_h\subseteq\Pi$ is a closed convex set, 
$h(\,\cdot\,)$ is strictly convex, twice differentiable in $\interior(\Pi_h)$, and is such that $\lim_{\pi\to\partial\Pi_h}||\tilde\nabla h(\pi)||=+\infty$, 
where $\tilde\nabla:=(I-|\mc P|^{-1}\onebf\onebf')\nabla$ is the projected gradient on $\mc S(\mc P)$\footnote{We shall use the notation $\Phi:=I-|\mc P|^{-1}\onebf\onebf'\in\R^{\mc P\times\mc P}$ to denote the corresponding projection matrix.}.
\end{assumption} 

As a result, the aggregate path preference $\pi(t)$ evolves as 
\be\label{slowscaleeqn}\frac{\de}{\de t}\pi=\eta\l(F^h(f)-\pi\r)\,.\ee
The perturbed best response function $F^h(f)$ provides an idealized description of the behavior of drivers whose decisions are based on inexact information about the state of the network. 
In particular, it can be shown that the form of $F^h(f)$ given in (\ref{Fhdef}) is equivalent to the minimization over paths $p\in\mc P$ of the expected delay $\sum_{e \in \mc E} A_{ep}T_e(f_e)$ corrupted by some (admissible) stochastic perturbation (see e.g.~\cite{Hofbauer.Sandholm:07}). 

It is easy to establish that the perturbed best response $F^h(f)$ is  continuously differentiable on $\mc F$. Moreover, it is well known \cite{Sandholm:10} that, as $\|h\|_\infty \downarrow 0$, and $\Pi_h\uparrow\ov\Pi$, the perturbed best response $F^h(f)$ converges to the set $\argmin\{\omega'A'T(f):\,\omega\in \Pi\}$ of best responses.\footnote{Here, the convergences $\Pi_h\uparrow\ov\Pi$, and $\{F^h(f)\}\to\argmin\{\omega'A'T(f):\,\omega\in \Pi\}$ are intended to hold in the Hausdorff metric. (see, e.g., \cite[Def.~4.4.11]{AmbrosioTilli})}
\begin{example}
Assume that $C_e>1$ for all $e\in\mc E$. Then, an example of perturbed best response satisfying Assumption \ref{assumption3} is the logit function with noise level $\beta>0$, which is defined as
\begin{equation}
\label{eq:logit}
F^h_p(f)=\frac{\exp(-\beta (A'T(f))_p)}{\sum_{q \in \paths}\exp(-\beta (A'T(f))_q)}\,,\qquad p\in\mc P\,.
\end{equation}
This corresponds to the admissible perturbation function  $h(\omega)=-{\beta}^{-1}\sum_p\omega_p\log\omega_p$. 
For any fixed $f\in\mc F$, one has that $\lim_{\beta\to+\infty}F^h(f)$, with $F^h(f)$ as defined in \eqref{eq:logit}, is a uniform distribution over the set $\argmin\{(A'T(f))_p:\,p\in\mc P\}$.  
We refer the reader to \cite{Fudenberg.Levine:08} for more on the connection between $F^h$ characterized by Assumption \ref{assumption3} and smooth best response functions.
\end{example}

\begin{remark}\label{remarkPih}
In the evolutionary game theory literature, e.g., see \cite{HofbauerSigmund:03,Sandholm:10}, the domain of an admissible perturbation function $h$, as well as the one of the minimization in the right-hand side of (\ref{Fhdef}), is typically assumed to be the whole simplex $\mc S(\mc P)$, instead of a closed polytope $\Pi_h\subseteq\Pi\subseteq\mc S(\mc P)$. 
Notice that, as already observed in Sect.~\ref{sect:networkcharacteristics}, when $C_e>1$ for every $e\in\mc E$, $\Pi=\mc S(\mc P)$ is a closed polytope, so that one can choose $\Pi_h=\Pi$. Therefore, in this case, Assumption \ref{assumption3} does not introduce any additional restriction with respect to such theory. 

On the other hand, when $C_e\le1$ for some $e\in\mc E$, then the inclusions of $\Pi_h\subset\Pi\subset\mc S(\mc P)$ are both strict, so that Assumption \ref{assumption3} does introduce additional restrictions on the admissible perturbations. However, it is worth observing that, in a classic evolutionary game theoretic framework, the dynamics of the aggregate path preference would be autonomous rather than coupled to the one of the actual flow. In particular, perturbed best response dynamics in that framework would read as 
\be\label{standardEGTeqn}\frac{\de}{\de t}\pi=F^h(f^{\pi})-\pi\,,\ee rather than as in (\ref{slowscaleeqn}). For such dynamics, the fact that $T_e((A\pi)_e)=+\infty$ whenever $(A\pi)_e\ge C_e$, can be shown to imply that $\pi(t)$ would reach a compact $\Pi_h\subseteq\Pi$ in some finite time and never leave it. In contrast, in the two time-scale model of coupled dynamics considered in this paper (see (\ref{coupleddyn})), such more restrictive assumption is needed in order to ensure the same property for the trajectories of $\pi(t)$ (see Lemma \ref{lemmaawayC1}).
\end{remark}

\subsubsection*{Local route decisions}
We now describe the local route decisions, characterizing the fraction of drivers choosing each link $e \in \edgeset_{v}^+$ when traversing a non-destination node $v$.
Such a fraction will be assumed to be a continuously differentiable function $G^v_e(f_{\mc E^+_v},\pi)$ of the local traffic flow $f_{\mc E^+_v}:=\{f_e:\,e\in\edgeset^+_v\}$, as well as of the current aggregate path preference $\pi$.
We shall refer to 
\be\label{Gdef}G^v:\mc F_v\times\Pi\to\mc S(\mc E^+_v)\ee
as the \emph{local decision function} at node $v\in\{0,1,\ldots, n-1\}$, 
and assume that it satisfies the following: 
\begin{assumption}\label{assumption4} For all $0\le v<n$, and $\pi\in\Pi$, 
$$\l(\sum\nolimits_{j\in\mc E^+_v}f_j^{\pi}\r)G^v_e\left(f^{\pi}_{\mc E^+_v},\pi \right)=f_e^{\pi}\,,\qquad\forall e\in\edgeset_v^+\,.$$
\end{assumption}
\begin{assumption}\label{assumption5} For all $0\le v<n$, $\pi\in\Pi$, and $f_{\mc E^+_v}\in\mc F_v$,  
$$\frac{\partial}{\partial f_e }G_j^v(f_{\mc E^+_v},\pi)\geq 0\,,\qquad \forall j\ne e\in\edgeset^+_v\,.$$ 
\end{assumption}

Assumption \ref{assumption4} is a consistency assumption. It postulates that, when the locally observed flow coincides with the one associated to the aggregate path preference $\pi$, drivers choose to join link $e\in\mc E^+_v$ with frequency equal to the ratio between the flow $f^{\pi}_e$ and the total outgoing flow $\sum_{j\in\mc E^+_v}f^{\pi}_j$.

Assumption \ref{assumption5} instead models the drivers' myopic behavior in response to variations of the local congestion levels. It postulates that, if the congestion on one link increases while the congestion on the other links outgoing from the same node is kept constant, the frequency with which each of the other outgoing links is chosen does not decrease. It is worth observing that Assumption \ref{assumption5} is reminiscent of Hirsch's notion of \emph{cooperative dynamical system} \cite{Hirsch:82,Hirsch:85}. 

\begin{example}\label{ex:i-logit}An example of local decision function $G^v$ satisfying Assumptions \ref{assumption4} and \ref{assumption5} is the i-logit function. The i-logit route choice with sensitivity $\gamma>0$ is given by \begin{equation}\label{eq:i-logit}G^v_e(f_{\mc E^+_v},\pi)=\frac{f^{\pi}_e\exp(-\gamma(f_e-f_e^{\pi}))}{\sum_{j\in\mc E_v^+} f^{\pi}_j \exp(-\gamma(f_j-f_j^{\pi}))}\,,\end{equation} for every $e\in\mc E_v^+$, $0\le v<n$. 
\end{example}

For every non-destination node $v\in\{0,1,\ldots,n-1\}$, and outgoing link $e\in\mc E^+_v$, conservation of mass implies that 
\be\label{fastscaledyn}\frac{\de}{\de t}\rho_e=H_e(f,\pi)\,,\qquad 
H_e(f,\pi):=\l\{\ba{lcl}G^v_e(f_{\mc E^+_v},\pi)-f_e&\se&v=0\\(\sum_{j\in\mc E^-_v}f_j)G^v_e(f_{\mc E^+_v},\pi)-f_e&\se&1\le v<n\,.\ea\r.\ee

\subsection{Objective of the paper and main result}
The objective of this paper is to study the evolution of the coupled dynamics 
\be\label{coupleddyn}
\l\{
\ba{l}\ds\frac{\de}{\de t}\pi=\eta\l(F^h(f)-\pi\r)\\[10pt]
\ds\frac{\de}{\de t}\rho=H(f,\pi)\,,
\ea\r.
\ee
where $F^h$ is the perturbed best response function defined in (\ref{Fhdef}), $\eta>0$ is the rate at which global information becomes available, $H(f,\pi)=\{H_e(f,\pi):\,e\in\mc E\}$, with $H_e$ defined in (\ref{fastscaledyn}), and $f$ and $\rho$ are related by the functional dependence (\ref{flowdensityeqn}). 
In particular, our analysis will focus on the double limiting case of small $\eta$ and small $h$.
We shall prove that, in such limiting regime, the long-time behavior of the system is approximately at Wardrop equilibrium \cite{Wardrop:52,Patriksson:94}. The latter is a configuration in which the delay is the same on all the paths chosen by a nonzero fraction of the drivers.  More formally, one has the following: 
\begin{definition}[Wardrop Equilibrium]
An admissible flow vector $f^{W}\in\mc F$ is a \emph{Wardrop equilibrium} if $f^{W}=A\pi$ for some $\pi\in\Pi$ such that, for all $p\in\mc P$, 
\be\label{WEeqn}\pi_p>0\qquad\Longrightarrow\qquad\l(A'T(A\pi)\r)_p\le\l(A'T(A\pi)\r)_q\,,\qquad \forall q\in\mc P\,.\ee
\end{definition}

\noindent Existence and uniqueness of a Wardrop equilibrium are guaranteed by the following standard result:
\begin{proposition}[Existence and uniqueness of Wardrop equilibrium]
\label{prop:unique}
Let Assumptions \ref{assumption1}-\ref{assumption2.5} be satisfied. Then, there exists a unique Wardrop equilibrium $f^{W}\in\mc F$.
\end{proposition}
\ksmargin{proof needs more details}
\proof
It follows from Assumption \ref{assumption2} that, for every $e \in \linkset$, the delay function $\delayfunc_e(f_e)$ is  continuous, strictly increasing, and such that $T_e(0)>0$. The proposition then follows by applying Theorems 2.4 and 2.5 from \cite{Patriksson:94}.
\qed  
 
The following is the main result of this paper. 
It will be proved in Sect.~\ref{sec:proofs} using a singular perturbation approach.
\begin{theorem}
\label{thm:main-unperturbed}
Let Assumptions \ref{assumption1}--\ref{assumption5} be satisfied. 
Then, for every initial condition $\pi(0)\in\interior(\mc S(\mc P))$, $\rho(0)\in(0,+\infty)^{\mc E}$, there exists a  unique solution of (\ref{coupleddyn}). Moreover, there exists a perturbed equilibrium flow $f^h \in \mc F$ such that,  for all $\eta>0$, 
\be\label{limsup<=eta}
\limsup_{t\to+\infty}||f(t)-f^h||\le\delta(\eta)\,,
\ee
where $\delta(\eta)$ is a nonnegative-real-valued, nondecreasing function of $\eta>0$, such that $\lim_{\eta\downarrow0}\delta(\eta)=0$. 
Moreover, for every sequence of admissible perturbations $\{h_k:\,k\in\N\}$ such that 
$\lim_{k t\ +\infty} ||h_k||_{\infty}=0$, and $\lim_{k \to +\infty}\Pi_{h_k}=\ov{\Pi}$, one has
\be\label{limh0}\lim_{k\to+\infty}f^{h_k}=f^{W}\,.\ee
\end{theorem}

Theorem \ref{thm:main-unperturbed} states that, in large time limit, the flow vector $f(t)$ approaches a neighborhood of the Wardrop equilibrium, whose size vanishes as both the time-scale ratio $\eta$ and the perturbation norm $||h||_{\infty}$ vanish. While a qualitatively similar result is known to hold \cite{Sandholm:10} in a classic evolutionary game theoretic framework (i.e., neglecting the traffic dynamics, and assuming it is instantaneously equilibrated, as in the ODE system (\ref{standardEGTeqn})), the significance of the above is to show that an approximate Wardrop equilibrium configuration is expected to emerge also in our more realistic model of two-time scale dynamics. Therefore, our results provide a stronger evidence in support of the significance of Wardrop's postulate of equilibrium for a transportation network. In fact, they may be read as a sort of robustness of such equilibrium notion with respect to non-persistent perturbations.  

\section{Proofs}
\label{sec:proofs}
In this section, Theorem \ref{thm:main-unperturbed} is proved. First, observe that, thanks to the continuous differentiability of $F^h$, $G^v$, and $\mu$, standard analytical arguments imply the existence and uniqueness of a solution of the the initial value problem associated to the system (\ref{coupleddyn}), with initial condition $\rho(0)\in(0,+\infty)^{\mc E}$, $\pi(0)\in\interior(\Pi)$. 

In order to prove the rest of the statement, we shall adopt a singular perturbation approach (e.g., see \cite{Khalil:96}), viewing the traffic density $\rho$ (or, equivalently, the traffic flow $f$) as a fast transient, and the aggregate path preference $\pi$ as a slow component. Hence, we shall first think of $\pi$ as quasi-static (i.e., `almost a constant') while analyzing the fast-scale dynamics \eqref{fastscaledyn}, and then assume that $f$ is `almost equilibrated', i.e.~close to $f^{\pi}$, and study the slow-scale dynamics \eqref{slowscaleeqn} as a perturbation of (\ref{standardEGTeqn}). We shall proceed by proving a series of intermediate technical results, gathered in the following subsections.

Before proceeding, we introduce some notation to be used throughout the section. 
Let $$\rho_e^{\pi}:=\mu_e^{-1}(f_e^{\pi})\,,\qquad\sigma_e:=\sgn\l(\rho_e-\rho_e^{\pi}\r)=\sgn\l(f_e-f_e^{\pi}\r)$$
denote, respectively, the density corresponding to the flow associated to the path preference $\pi$, and the sign of the difference between it and the actual density $\rho_e$. 
Finally, fix some $\alpha\in(0,1)$, and define
$$V(f,\pi):=\sum_{v=0}^{n-1}\alpha^v\sum_{e\in\mc E^+_v}|f_e-f^{\pi}_e|\,,\qquad\qquad 
W(\rho,\pi):=\sum_{v=0}^{n-1}\alpha^v\sum_{e\in\mc E^+_v}|\rho_e-\rho^{\pi}_e|.$$

\subsection{Stability of the fast-scale dynamics}
We gather here a few properties of the fast-scale dynamics. 
Our results will essentially amount to showing that $V(f,\pi)$ and $W(\rho,\pi)$ are Lyapunov functions for the fast-scale dynamics (\ref{fastscaledyn}) with stationary path preference $\pi$. 

The following result is a consequence of Assumptions \ref{assumption4} and \ref{assumption5} on the drivers' local decision function. 

\begin{lemma}\label{lemma0}
For all $\pi\in\Pi$, $v \in \{0,\ldots,n-1\}$, and $f_{\mc E^+_v}\in\mc F_v$,
$$\sum_{e\in\mc E^+_v}\sigma_e\l(\lambda^{\pi}_v G^v_e(f_{\mc E^+_v},\pi)-f_e^{\pi}\r)\le0\,,$$
where $\lambda^{\pi}_v:=\sum_{e\in\mc E^+_v}f_e^{\pi}$. 
\end{lemma}
\proof 
Throughout this proof, the explicit dependence of $G^v_e$ on $\pi$ will be dropped. 
Define $\mc J:=\{e\in\mc E^+_v:\,f_e>f_e^{\pi}\}$, $\mc K:=\{e\in\mc E^+_v:\,f_e<f_e^{\pi}\}$, and let $G_{\mc J}:=\sum_{j\in\mc J} G_j^{v}$, $G_{\mc K}:=\sum_{k\in\mc K} G_k^{v}$, and $G_{\mc J^c}:=\sum_{e\in\mc E^+_v\setminus\mc J} G_e^{v}$. 
First, observe that, since $\sum_{e \in \edgeset_v^+} G_e^v=1$, one has that $\nabla G_{\mc J}=-\nabla G_{\mc J^c}$. 
Now, we are going to show that  
\be\label{in0} G_{\mc J}(f^{\pi}_{\mc E^+_v}) - G_{\mc J}(f_{\mc E^+_v}) \ge0\,,\ee
by writing the difference above as a path integral of $\nabla G_{\mc J}(\,\cdot\,)$ first along the segment $S_{\mc J}$ from $f_{\mc E^+_v}$ to the point $f^*\in\R_+^{\mc E^+_v}$ with $f_j^*:=f_j^{\pi}$, for $j\in\mc J$ and $f_e^*:=f_e$ for $e\in\mc E^+_v\setminus\mc J$, and then along the segment $S_{\mc K}$ from $f^*$ to $f^{\pi}$. In this way, one gets:
\begin{equation}
\label{eq:path-integral}
\ba{rcl}G_{\mc J}(f^{\pi}_{\mc E^+_v}) - G_{\mc J}(f_{\mc E^+_v})&=& 
\ds\int_{S_{\mc J}}\nabla G_{\mc J}( \tilde f_{\mc E^+_v})\cdot\de\tilde f_{\mc E^+_v}+ 
\int_{S_{\mc K}}\nabla G_{\mc J}(\tilde f_{\mc E^+_v})\cdot\de\tilde f_{\mc E^+_v} \\[7pt]
&=&\ds
-\int_{S_{\mc J}}\nabla G_{\mc J^c}(\tilde f_{\mc E^+_v})\cdot\de\tilde f_{\mc E^+_v}+ 
\int_{S_{\mc K}}\nabla G_{\mc J}(\tilde f_{\mc E^+_v})\cdot\de\tilde f_{\mc E^+_v} 
\,. \ea\end{equation}
Assumption \ref{assumption5} implies that $\partial G_{\mc J^c}/\partial\rho_j\ge0$ for all $j\in\mc J$, and $\partial G_{\mc J}/\partial\rho_k\ge0$ for all $k\in\mc K$. 
In turn, this implies that $\nabla G_{\mc J^c}\cdot\de\tilde f_{\mc E^+_v}\le0$ along $S_{\mc J}$, and $\nabla G_{\mc J}\cdot\de\tilde f_{\mc E^+_v}\ge0$ along $S_{\mc K}$. 
This and (\ref{eq:path-integral}) prove (\ref{in0}). In a very similar fashion, one proves that 
\be\label{in1}G_{\mc K}(f_{\mc E^+_v}) - G_{\mc K}(f_{\mc E^+_v}^{\pi})\ge0\,.\ee
Now, observe that Assumption \ref{assumption4} implies that $\lambda_v^{\pi}G^v_e(f^{\pi}_{\mc E^+_v},\pi)=f_e^{\pi}$. 
From this,  (\ref{in0}), and (\ref{in1}), it follows that 
$$\ba{rcl}
0&\ge&\ds
\lambda_v^{\pi}\l(G_{\mc J}(f_{\mc E^+_v})-G_{\mc J}(f_{\mc E^+_v}^{\pi})\r)
-\lambda_v^{\pi}\l(G_{\mc K}(f_{\mc E^+_v})-G_{\mc K}(f_{\mc E^+_v}^{\pi})\r)\\[7pt]
&=&\ds
\sum_{e\in\mc E^+_v}\sigma_e\l(\lambda_v^{\pi} G^v_e(f_{\mc E^+_v})-\lambda_v^{\pi} G^v_e(f_{\mc E^+_v}^{\pi})\r)\\[7pt]
&=&\ds
\sum_{e\in\mc E^+_v}\sigma_e\l(\lambda_v^{\pi} G^v_e(f_{\mc E^+_v})-f_e^{\pi}\r)
\,,\ea$$
which proves the claim. 
\qed

We now proceed to analyzing, for a fixed global decision $\pi \in \simplex$, the fast scale dynamics (\ref{fastscaledyn}). 
Let 
$$V^+_v(f,\pi):=\sum\nolimits_{e\in\mc E^+_v}\l|f_e^{\pi}-f_e\r|\,,\qquad v=0,1,\ldots,n-1\,,$$
be the $l_1$-distance between the current flows on the outgoing links of $v$, and the flow associated to the aggregate path preference $\pi$, and
$$V^{-}_v(f,\pi):=\l|\lambda_v^{\pi}-\lambda_v^-\r|\,,\qquad v=1,2,\ldots,n\,,$$
with $\lambda^{\pi}_v:=\sum_{e\in\mc E^+_v}f_e^{\pi}$ and $\lambda_v^-:=\sum_{e \in \mc E_v^-} f_e$,
be the absolute difference between the current flow incoming in node $v$, and the one associated to the aggregate path preference $\pi$. Also, let $V^-_0(f,\pi):=0$. 


\begin{lemma}\label{lemmalocal}
For all  $v=0,1,\ldots,n-1$, $\pi\in\Pi$, and $f\in\mc F$, 
$$\sum_{e\in\mc E^+_v}\sigma_eH_e(f,\pi)\le-V^+_v(f,\pi)+V^-_v(f,\pi)\,.$$
\end{lemma}
\proof
Writing $G^v_e$ for $G^v_e(f_{\mc E^+_v},\pi)$, and using Lemma \ref{lemma0}, one gets  that 
$$\ba{rcl}\summ_{e\in\mc E^+_v}\!\!\sigma_eH_e(f,\pi)\!\!
&=&\!\!\!\!\summ_{e\in\mc E^+_v}\sigma_e\l(\lambda^-_v G^v_e-f_e\r)\\[10pt] 
&=&\!\!\!\!\summ_{e\in\mc E^+_v}\!\!\sigma_e(\lambda_v^--\lambda_v^{\pi})G^v_e
+\!\!\summ_{e\in\mc E^+_v}\!\!\sigma_e\l(\lambda_v^{\pi}G^v_e-f_e^{\pi}\r)
+\!\!\summ_{e\in\mc E^+_v}\!\!\sigma_e\l(f_e^{\pi}-f_e\r)\\[10pt]
&\le&|\lambda_v^--\lambda^{\pi}_v|-\summ\nolimits_{e\in\mc E^+_v}\l|f_e^{\pi}-f_e\r|\\[10pt]
&=&\!\!-V_v^+(f,\pi)+V^-_v(f,\pi)\,,
\ea$$
which proves the claim. 
\qed

By combining Lemma \ref{lemmalocal}, and Assumption \ref{assumption1}, one gets the result below. Recall that we are using the convention ${\de}|x|/{\de x}=\sgn(x)$, for all $x\in\R$. 
\begin{lemma}\label{lemmaexpstabslowscale}
For every $f=\mu(\rho)\in\mc F$, and $\pi\in\Pi$, 
$$\nabla_{\rho}W(\rho,\pi)'H(f,\pi)\le-(1-\alpha)V(f,\pi)\,.$$
\end{lemma}
\proof
Observe that, thanks to the acyclicity of the graph as per Assumption~\ref{assumption1}, if $e\in\mc E^-_v\cap\mc E^+_w$ for some nodes $v$ and $w$, then necessarily $v\ge w+1$. Since $\alpha<1$, it follows that $$\alpha^v\1_{\mc E_v^-}(e)\1_{\mc E_w^+}(e)\le \alpha^{w+1}\1_{\mc E_v^-}(e)\1_{\mc E_w^+}(e)\,,$$
for every  $1\le v\le n$, and $0\le w\le n-1$. Hence, 
$$
\ba{rcl}\ds\summ_{0\le v<n}\alpha^vV^-_v(f,\pi)
&\le &
\ds\summ_{0\le v<n}\summ_{e\in\mc E_v^-}\alpha^v\l|f_e-f_e^{\pi}\r|\\[7pt]
&=&\ds\summ_{1\le v<n}\summ_{0\le w<n}\summ_{e\in\mc E}\alpha^v\1_{\mc E_v^-}(e)\1_{\mc E_w^+}(e)\l|f_e-f_e^{\pi}\r|\\[7pt]
&\le&\ds\summ_{0\le w<n}\alpha^{w+1}\summ_{e\in\mc E}\1_{\mc E_w^+}(e)\l|f_e-f_e^{\pi}\r|\summ_{1\le v<n}\1_{\mc E_v^-}(e)\\[7pt]
&\le&\alpha\ds\summ_{0\le w<n}\alpha^{w}\summ_{e\in\mc E^+_w}\l|f_e-f_e^{\pi}\r|\\[7pt]
&=&\alpha V(f,\pi)\,,
\ea
$$
where the last inequality follows from the fact that $\sum_{v=1}^{n-1}\1_{\mc E_v^-}(e)\le\sum_{v=1}^n\1_{\mc E_v^-}(e)=1$. 
Thus, Lemma \ref{lemmalocal} implies that 
$$
\ba{rcl}
\ds\nabla_{\rho}W(\rho,\pi)'H(f,\pi)
\!\!\!\!&=&\ds
\summ_{0\le v<n}\alpha^v\summ_{e\in\mc E^+_v}\sigma_e H_e(f,\pi)\\[7pt]
&\le&
\ds\summ_{0\le v<n}\alpha^vV^-_v(f,\pi)-\summ_{0\le v<n}\alpha^vV_v^+(f,\pi)\\[7pt]
&\le&\alpha V(f,\pi)-V(f,\pi)\,,
\ea
$$
which proves the claim.\qed

\subsection{Boundedness of the traffic densities} 
We shall now prove a couple of results guaranteeing that the traffic density on every link $e\in\mc E$ remains bounded in time. 
We start with the following result, guaranteeing that, on every link $e\in\mc E$, the flow associated to the current path preference, $f^{\pi}_e(t)$, stays eventually bounded away from the maximum flow capacity $C_e$. 
Its proof relies on Assumption \ref{assumption3}. Recall that our formulation allows for both the cases of finite and infinite maximum flow capacity on a link.

\begin{lemma}\label{lemmaawayC1}
For every admissible perturbation $h$, there exists $t_0\in\R_+$, and, for every link $e \in \mc E$, a positive finite constant $\ov C_e$, dependent on $h$ but not on $\eta$, such that, for every initial condition $\pi(0)\in \interior \left(\mc S(\mc P)\right)$, $\rho(0)\in(0,+\infty)^{\mc E}$,
$$f^{\pi}_e(t)\le \ov C_e < C_e \,,\qquad\forall t\ge t_0\,, \qquad \forall \, e \in \mc E.$$
\end{lemma}
\proof 
The fact that $f_e^\pi(t) \leq 1$ for all $e \in \mc E$ follows from the fact that the arrival rate at the origin is unitary. Therefore, for all $e \in \mc E$ with $C_e > 1$ or $C_e=\infty$ in particular, the lemma follows trivially with $\ov C_e=1$ and $t_0=0$. 
We now prove the lemma for all $e \in \mc E$ with $C_e < \infty$.
Recall that, by Assumption \ref{assumption3}, the domain of the admissible perturbation $h$ is a closed set $\Pi_h\subset\interior(\Pi)$. 
This, in particular implies that, for all $e \in \mc E$ with $C_e<\infty$, 
$$\boundconst_e:=C_e-\sup\{(A\tilde\pi)_e:\,\tilde\pi\in\Pi_h\}>0\,.$$ 
For every link $e\in\mc E$ with $C_e < \infty$, it follows from (\ref{Fhdef}) that 
\be\label{supFe}
\ba{rcl}
C_e-\boundconst_e&=&\sup\{(A\tilde\pi)_e:\,\tilde\pi\in\Pi_h\}\\[3pt]
&\ge&\sup\l\{\l(A\argmin\{\tilde\pi'A'T(f)+h(\tilde\pi):\,\tilde\pi\in\Pi_h\}\r)_e:\,f\in\mc F\r\}\\[3pt]
&=&\sup\l\{\l(AF^h(f)\r)_e:\,f\in\mc F\r\}\,.\ea\ee
Hence, one has for every link $e\in\mc E$ with $C_e < \infty$,
$$\frac{\de}{\de t}f^{\pi}_e(t)=\eta\l(A(F^h(f(t))-\pi(t))\r)_e\le\eta\l(C_e-\boundconst_e-f^{\pi}_e\r)\,.$$
Then, Gronwall's inequality implies that 
$$f_e^{\pi}(t)\le f_e^{\pi}(0)e^{-\eta t}+(1-e^{-\eta t})(C_e-\boundconst_e)\le e^{-\eta t}+C_e-\boundconst_e\,,$$
for all $t\ge0$ and for every link $e\in\mc E$ with $C_e < \infty$.  The lemma for $e \in \mc E$ with $C_e < \infty$ now follows from the above, by choosing, e.g., $\ov C_e:=C_e - \kappa$ with $\boundconst:=\frac12\min\{\boundconst_e:\,e\in\mc E \text{ s.t. } C_e<\infty\}$, and $t_0:=-\eta^{-1}\log\boundconst$.
\qed

The following result shows that the actual flow $f_e(t)$ also stays bounded away from the maximum flow capacity $C_e$.
\begin{lemma}\label{lemmaawayC2}
For every admissible perturbation $h$, there exists positive finite constants $\eta^*$, and $\tilde C_e$, for every $e \in \mc E$, such that, for every $\eta<\eta^*$, and every initial condition $\pi(0)\in \interior \left(\mc S(\mc P)\right)$, $\rho(0)\in(0,+\infty)^{\mc E}$,
$$f_e(t)\le \tilde C_e < C_e \,,\qquad \forall t\ge0,\, \quad \forall \, e \in \mc E.$$
\end{lemma}
\proof
Let $\zeta(t):=W(\rho(t),\pi(t))$, and $\chi(t):=V(f(t),\pi(t))$. 
Observe that, thanks to Lemma \ref{lemmaawayC1}, there exists a positive real constant $\ov C_e$ for every $e \in \mc E$, and $t_0\ge0$, such that, for every $t\ge t_0$, 
\be\label{rhopifinite}\rho^{\pi}_e(t)\le\rho_e^*\,,\qquad \rho_e^*:=\mu_e^{-1}(\ov C_e)\,,\qquad\forall e\in\mc E\,.\ee 
Since $\rho^{\pi}_e(t)\ge0$, the above implies that, if $|\rho_e(t)-\rho^{\pi}_e(t)|\ge2\rho_e^*$ for some $e \in \mc E$ and $t\ge t_0$, then necessarily for that $e \in \mc E$, we have that $f_e(t)-f_e^{\pi}(t)\ge\chi_e^*$, where $\chi_e^*:=\mu_e(2\rho_e^*)-\ov C_e>0$. Now, let $\zeta^*:=2 |\mc E|\max\{\rho^*_e:e\in\mc E\}$, and $\chi^*:=\alpha^{n-1}\min\{\chi^*_e:e\in\mc E\}$. Note that $W(\rho,\pi)\le|\mc E|\max\{|\rho_e-\rho^{\pi}_e)|:e\in\mc E\}$, and $V(f,\pi)\ge\alpha^{n-1} |f_e-f_e^{\pi}|$ for every $e \in \mc E$. Therefore, it follows that, for any $t \geq t_0$, if $\zeta(t) \geq \zeta^*$, then for some $e' \in \mc E$, we have that $|\rho_{e'}(t)-\rho_{e'}^\pi| \ge 2 \rho^*_{e'}$ for $t \ge t_0$. This in turn implies that $\chi(t) \geq \chi_{e'}^* \geq \chi^*$. Therefore, in summary,
\be\label{implication}\zeta(t)\ge\zeta^*\quad\Longrightarrow\quad\chi(t)\ge\chi^*\,,\qquad\forall t\ge t_0.\ee 

On the other hand, observe that (\ref{rhopifinite}) implies that there exists some $\ell>0$ such that 
$$\sum_{0\le v<n}\alpha^v\sum_{e\in\mc E^+_v}\frac1{\mu_e'(\rho^{\pi}_e(t))}\le \ell\,,\qquad\forall t\ge t_0\,.$$
By combining the above with Lemma \ref{lemmaexpstabslowscale}, one finds that, for any $u,t \geq t_0$, 
\be\label{gammat}
\ba{rcl}\zeta(t)-\zeta(u)&=&
\ds\int_{u}^t\sum_{0\le v<n}\alpha^v\sum_{e\in\mc E^+_v}\sigma_e\l(\frac{\de}{\de s}\rho_e-\frac{\de}{\de s}\rho^{\pi}_e\r) \de s\\
&\le&\ds\int_{u}^t\nabla_{\rho}W(\rho,\pi)'H(f,\pi) \de s\\
&&+\ds\int_{u}^t\sum_{0\le v<n}\alpha^v\sum_{e\in\mc E^+_v}\frac{\eta}{\mu_e'(\rho^{\pi}_e)}\l|(AF^h(f^\pi))_e-(A\pi)_e\r| \de s\\
&\le&\ds\int_{u}^t\l(-(1-\alpha)\chi(s)+2 \eta \ell\r)\de s\,.
\ea\ee

Now, define $\eta^*:=(1-\alpha)\chi^*/(2\ell)$. Assume, by contradiction, that $\limsup_{t \to \infty}$ $f_e(t) \ge C_e$ for some $e\in\mc E$. Since $f_e(t)=\mu_e(\rho_e(t))<C_e$ for every $t\ge0$, this implies that $\limsup_{t\to+\infty}\rho_e(t)=+\infty$. 
This, together with (\ref{rhopifinite}) implies that $\limsup_{t\to+\infty}$ $\zeta(t)=+\infty$.
Then, in particular, the set $\mc T:=\{t>0:\zeta(t)>\zeta(s)\,,\forall s<t\}$ is an unbounded union of open intervals, with
$\lim_{{t\in\mc T,t\to+\infty}}\zeta(t)=+\infty$.
This, and (\ref{implication}) imply that there exists a non negative constant $t^* \geq t_0$ such that 
$$\chi(t)\ge\chi^*\,,\qquad \forall t\ge t^*\,.$$
For every $\eta<\eta^*$, Equation (\ref{gammat}) and the above give
$$\zeta(t)-\zeta(u)=\ds\int_{u}^t\l(-(1-\alpha)\chi(s)+2 \eta \ell\r)\de s
\le \ds\int_{u}^t\l(-(1-\alpha)\chi^*+2 \eta \ell\r)\de s<0$$
for every $t>u\ge t^*$ such that $t$ and $u$ belong to the same connected component of $\mc T$.
But this contradicts the definition of the set $\mc T$.  
Hence, if $\eta<\eta^*$, then $\limsup_{t\to+\infty}f_e(t)<C_e$ for every $e\in\mc E$. Since on every compact time interval $\mc I\subseteq[0,+\infty)$, one has $\sup_{t\in\mc I}f_e(t)=f_e(\hat t)<C_e$ for some $\hat t\in\mc I$, the foregoing implies the lemma.
\qed

The result below is a consequence of Lemma \ref{lemmaawayC2}, and will prove useful in the sequel. 
\begin{proposition}\label{prop} 
There exists $K>0$, and $t_1 \ge0$ such that, for every initial condition $\pi(0)\in \interior \left(\mc S(\mc P)\right)$, $\rho(0)\in(0,+\infty)^{\mc E}$, $||\tilde\nabla_{\pi}h(\pi(t))||\le K$ for all $t\ge t_1$.
\end{proposition}
\proof
First, observe that, thanks to Lemma \ref{lemmaawayC2}, there exists $T^*>0$ such that $||T(f(t))||\le T^*$ for all $t\ge0$. Thanks to this, and Assumption \ref{assumption3}, one has that $F^h(f(t))\in\interior(\Pi_h)$, and $\tilde\nabla_{\pi}h(F^h(f(t)))=-\Phi A'T(f(t))$, where recall that $\Phi=I-|\mc P|^{-1}\onebf\onebf'$ is the projection matrix corresponding to the projected gradient with respect to $\pi$ on $\mc S(\mc P)$. Hence, $||\tilde\nabla_{\pi}h(F^h(f(t)))||\le ||\Phi|| ||A'||T^*$, which implies that there exists a convex compact $\mc K\subseteq\interior(\Pi_h)$ such that $F^h(f(t))\in\mc K$ for all $t\ge0$. Define $$\Delta(t):=\eta \l(1-e^{- \eta t}\r)^{-1}\int_0^te^{-\eta (t-s)}F^h(f(s))\de s\,.$$ As $\Delta(t)$ is an average of elements of the convex set $\mc K$, necessarily $\Delta(t)\in\mc K$ for all $t \geq 0$. Then, $\pi(t)=e^{-\eta t}\pi(0)+(1-e^{- \eta t})\Delta(t)$ approaches $\mc K$, which implies that, for large enough $t$, $\pi(t)\in \mc K_1 \subset  \interior(\Pi_h)$, where $\mc K_1$ is a closed subset of $\interior(\Pi_h)$ that contains $\mc K$.
Hence, after large enough $t$, say $t_1$, $\tilde\nabla_{\pi}h(\pi(t))$ stays bounded. 
\qed

\subsection{Estimating the distance between the current density and the one associated to the current path preference}
We analyze here the behavior in time of $W(\rho(t),\pi(t))$. 
First, we have the following result, characterizing the variation of $W(\rho,\pi)$ as a function of $\pi$. Recall that $\tilde\nabla_{\pi}=(I-|\mc P|^{-1}\onebf\onebf')\nabla_{\pi}$ denotes the projected gradient with respect to $\pi$ on $\mc S(\mc P)$. 

\begin{lemma}\label{lemmaexpstabbis} There exists $l>0$, and $t_0\ge0$, such that, for every initial condition $\pi(0)\in \interior \left(\mc S(\mc P)\right)$, $\rho(0)\in(0,+\infty)^{\mc E}$,
$$\tilde\nabla_{\pi}W(\rho(t),\pi(t))'(F^h(f(t))-\pi(t))\le\frac{2l}{1-\alpha}\,,\qquad\forall t\ge t_0\,.$$
\end{lemma}
\proof First, observe that, thanks to Lemma \ref{lemmaawayC1}, one has that there exists $t_0\ge0$ such that $l_e:=\sup\{1/\mu_e'(\rho_e^{\pi}(t)):\,t\ge t_0\}<+\infty$. Put $l:=\max\{l_e:\,e\in\mc E\}$. Then, for every path $p\in\mc P$, and every $t\ge t_0$, one has  
\be\label{ineq1}
\ba{rcl}\ds\l|\frac{\partial W(\rho,\pi)}{\partial\pi_p}\r|
&=&\Bigg|-\ds\sum_{0\le v<n}\alpha^v\summ_{e\in\mc E^+_v}\sigma_e\frac{\partial}{\partial_{\pi_p}}\rho_e^{\pi}\Bigg|\\
&=&\Bigg|\ds\sum_{0\le v<n}\alpha^v\summ_{e\in\mc E^+_v}\sigma_e\frac{\partial}{\partial_{\pi_p}}\mu_e^{-1}\l(\sum\nolimits_pA_{ep}\pi_p\r)\Bigg|\\
&\le&\ds\sum_{0\le v<n}\alpha^v\summ_{e\in\mc E^+_v}A_{ep}\frac1{\mu_e'(\rho^{\pi}_e)}\\
&\le&\ds\sum_{0\le v<n}\alpha^v\summ_{e\in\mc E^+_v}A_{ep}l_e\\
&\le&\ds\frac{l}{1-\alpha}\,,\ea
\ee
where the third inequality follows from the fact that, thanks to Assumption \ref{assumption1} on the acyclicity of the network, each path $p\in\mc P$ passes through at most one link $e\in\mc E^+_v$. 
Therefore, 
$$\ba{rcl}\ds\frac{2l}{1-\alpha}
&\ge&\ds\sum_pF^h_p(f)\l|\frac{\partial}{\partial\pi_p}W(\rho,\pi)\r|+\sum_p\pi_p\l|\frac{\partial}{\partial\pi_p}W(\rho,\pi)\r|\\[7pt]
&\ge&\ds\sum_pF^h_p(f)\frac{\partial}{\partial\pi_p}W(\rho,\pi)-\sum_p\pi_p\frac{\partial}{\partial\pi_p}W(\rho,\pi)\\[7pt]
&=&\ds\tilde\nabla_{\pi}W(\rho,\pi)'(F^h(f)-\pi)\,,\ea$$
where the first inequality follows upon recalling that both $F^h(f)$, and $\pi$ are probability vectors over the path set $\mc P$, and by using \eqref{ineq1}.
\qed

We can now combine Lemmas \ref{lemmaexpstabslowscale} and \ref{lemmaexpstabbis}, in order to get the following estimate of the behavior in time of $W(\rho(t),\pi(t))$. 
\begin{lemma}\label{lemmaperturb}
There exist $l>0$, $L>0$, $\eta^*>0$ and $t_0\ge 0$ such that, for every initial condition $\pi(0)\in \interior \left(\mc S(\mc P)\right)$, $\rho(0)\in(0,+\infty)^{\mc E}$, 
\begin{multline*}
W(\rho(t),\pi(t)) \\ \le\frac{2\eta l L}{(1-\alpha)^2}+\l(W(\rho(t_0),\pi(t_0))-\frac{2\eta l L}{(1-\alpha)^2}\r)e^{-\frac{1-\alpha}{L}(t-t_0)}\, \quad \forall t \geq t_0, \quad \forall \, \eta > \eta^*.
\end{multline*}
\end{lemma}
\proof
Define $x(t):=W(\rho(t),\pi(t))$. Notice that, thanks to Lemmas \ref{lemmaawayC1} and \ref{lemmaawayC2}, there exist $L>0$, $\eta^* > 0$ and $t_0\ge0$, such that, for any $\eta < \eta^*$, $|\rho_e(t)-\rho_e^{\pi}(t)|\le L |f_e(t)-f_e^{\pi}(t)|$ for every $e\in\mc E$, $t\ge t_0$. This in particular implies that $V(f(t),\pi(t)) \ge W(\rho(t),\pi(t))/L= x(t)/L$, for all $\eta < \eta^*$ and $t\ge t_0$. Observe that $W(\rho,\pi)$ is a Lipschitz function of $\rho$ and $\pi$, while both $\rho(t)$ and $\pi(t)$ are Lipschitz on every compact time interval. Therefore, $x(t)$ is Lipschitz on every compact time interval, and thus differentiable for almost every $t\ge t_0$. For every $t$ at which $x(t)$ is differentiable, Lemmas \ref{lemmaexpstabslowscale} and \ref{lemmaexpstabbis} imply that 
$$
\ba{rcl}
\ds\frac{\de}{\de t}x(t)\!\!
&\le&
\ds\nabla_{\rho}W(\rho,\pi)'H(f,\pi)+\eta\tilde\nabla_{\pi}W(\rho,\pi)'(F^h(f)-\pi)\\[7pt]
&\le&\ds-(1-\alpha)V(f,\pi)+\frac{2\eta l}{1-\alpha}\\[7pt]
&\le&\ds-\frac{(1-\alpha)}{L}x(t)+\frac{2\eta l}{1-\alpha}\,. 
\ea
$$
Then, the claim follows from Gronwall's inequality. 
\qed

\subsection{Proof of Theorem \ref{thm:main-unperturbed}}
We now proceed to proving Theorem \ref{thm:main-unperturbed}.
Consider the following candidate Lyapunov function:
\begin{equation}
\label{eq:potential-def}
\psi^h(\pi) := \sum_{e \in \mc E} \int_0^{f_e^{\pi}} \delayfunc_e \left( s \right) ds + h(\pi)\,.
\end{equation}
Since $T_e(f_e)$ is increasing, one has that each term $\int_0^{f_e^{\pi}} \delayfunc_e \left( f_e \right)\de f_e $ is convex in $f_e^{\pi}$. Hence, the composition with the linear map $\pi\mapsto f^{\pi}_e=\sum_{p \in \mc P}A_{ep}\pi_p$ is convex in $\pi$. Since $h(\pi)$ is strictly convex by assumption, one gets that $\psi^h(\pi)$ is strictly convex as well. Therefore, $\psi^h(\pi)$ admits a unique minimizer \be\label{pihdef}\pi^h:=\argmin\l\{\psi^h(\pi):\,\pi\in\Pi_h\r\}\,.\ee
Let  $f^h:=A\pi^h$. Then, we have the following: 
\begin{lemma}\label{lemmahto0}
Let $\{h_k:\,k\in\N\}$ be any sequence of admissible perturbation functions such that $\lim\limits_{k}||h_k||_{\infty}=0$, $\lim\limits_{k}\Pi_{h_k}=\ov\Pi$. Then, 
$$\lim_{k\to+{\infty}}f^{h_k}=f^{W}\,.$$
\end{lemma}
\proof
Write $\pi^k$ for $\pi^{h_k}$, $F^k$ for $F^{h_{k}}$, and $\Pi_k$ for $\Pi_{h_k}$.  Since $\{A\pi^k\}\subseteq A\Pi$, and $A\ov\Pi$ is compact, there exists a converging sub-sequence $\{A\pi^{k_j}:\,j\in\N\}$. Let us denote by $f^*:=\lim_{j}A\pi^{k_j}\in A\ov\Pi$ its limit, and choose some $\pi^*\in\ov\Pi$ such that $f^*=A\pi^*$. Notice that, since $\sup\{T_e(f^{\pi}_e):\,\pi\in\Pi_h\}<+\infty$, Assumption \ref{assumption3} implies that the minimizer in (\ref{pihdef}) has to be in the interior of $\Pi_h$. As a consequence, one finds that necessarily $\tilde\nabla_{\pi}h(\pi^{k_j})=-\Phi A'T(A\pi^{k_j})$, which in turn implies that
$F^{{k_j}}(A \pi^{k_j})=\pi^{k_j}$. Then, using (\ref{Fhdef}), one finds that 
\be\label{kj}(A\pi^{k_j})'T(A\pi^{k_j})+h_{k_j}(\pi^{k_j})\le(A\pi)'T(A\pi^{k_j})+h_{k_j}(\pi)\,,\ee
for all $\pi\in\Pi_{k_j}$. Now, fix any $\pi\in\Pi$. Since $\Pi_k\stackrel{k}{\to}\ov\Pi$, one has that $\pi\in\Pi_{k_j}$ for all sufficiently large values of $j$. Hence, passing to the limit as $j\to+\infty$ in (\ref{kj}), one finds that 
$$(\pi^*)'A'T(A\pi^*)\le\pi'A'T(A\pi^*)\,,\qquad\forall \pi\in\Pi\,.$$
In turn, the above can be easily shown to be equivalent to the condition (\ref{WEeqn}) characterizing Wardrop equilibria. From the uniqueness of the Wardrop equilibrium, it follows that necessarily $f^*=f^{W}$. Then the claim follows from the arbitrariness of the accumulation point $f^*$, and the compactness of $A\ov\Pi$. 
\qed

We shall now estimate the time derivative of $\psi^h(\pi)$ along trajectories of our dynamical system. For this, define $$\Gamma(t):=\psi^h(\pi(t)).$$ 
Then, one has 
\be\label{dgammadt0}
\ba{rcl}
\ds\frac{\de}{\de t}\Gamma(t)
&=&\ds\tilde\nabla_{\pi}\psi^h(\pi(t))'\frac{\de}{\de t}\pi\\[7pt]
&=&\ds\eta a'\l(F^h(f(t))-\pi(t)\r)\\[7pt]
&=&\ds\eta a'\l(F^h(A\pi(t))-\pi(t)\r)+\eta a'\l(F^h(f(t))-F^h(A\pi(t))\r)\,.
\ea
\ee
where  $a:=\Phi A'T(A\pi(t))+\tilde\nabla_{\pi}h(\pi(t))$. 

Lemma~\ref{lemmaperturb} implies that there exists $t_2 \geq 0$, $\eta^* > 0$ and $M_1>0$ such that, for any $\eta < \eta^*$, $W(\rho(t),\pi(t)) \leq \eta M_1$ for all $t \geq t_2$. From the definition of $W$, it also follows that $W(\rho,\pi) \geq \alpha^{n-1} \|\rho-\rho^\pi\|_1$ for all $\rho,\pi$. Moreover, following Assumption~\ref{assumption2}, with $\ov L := \max_{e \in \mc E} \frac{\de \mu_e}{\de \rho_e}(0)$, we also have that $\|f -A \pi \|_1 \leq \ov L \|\rho-\rho^\pi\|_1$ for all $f=\mu(\rho)$ and $\pi$. Combining all these relationships, one can see that there exists a $M>0$ such that, for any $\eta < \eta^*$,  
\begin{equation}
\label{eq:f-Api}
||f(t)-A\pi(t)||\le\eta M, \qquad \forall \, t \geq t_2.
\end{equation}
Moreover, recall that $F^h$ is differentiable on $\mc F$, and that, thanks to Lemmas \ref{lemmaawayC1} and \ref{lemmaawayC2}, for $\eta < \eta^*$, both $f(t)$ and $A\pi(t)$ are eventually confined in a compact $\mc K\subseteq\mc F$. This implies that 
$$||F^h(f(t))-F^h(A\pi(t))||\le K_1\eta$$ for some positive constant $K_1$, $\eta< \eta^*$ and sufficiently large values of $t$. 
On the other hand, Lemma \ref{lemmaawayC1} and Proposition \ref{prop} imply that both $T(A\pi(t))$ and $\tilde\nabla_{\pi}h(\pi(t))$ are eventually bounded, so that $||a||\le K_2$, for some positive constant $K_2$ and large enough $t$. 
It follows that the second addend in the last line of (\ref{dgammadt0}) can be bounded as 
\be\label{dgammadt1}\eta a'\l(F^h(f(t))-F^h(A\pi(t))\r)\le K\eta^2\,,\qquad \forall \eta< \eta^*, \quad \forall \, t\ge t_3 \,,\ee
for some sufficiently large but finite value of $t_3$, where $K=K_1K_2$. Now, observe that, for every $\pi$,  
$$\Phi A'T(A\pi)=-\tilde\nabla_{\pi}h\l(F^h(A\pi)\r)\,,$$
so that the first addend in the last line of (\ref{dgammadt0}) may be rewritten as 
\be\label{dgammadt2}
a'\l(F^h(A\pi(t))-\pi(t)\r)=-\Upsilon(\pi(t))\,,
\ee
where 
$$\Upsilon(\pi):=\l(\tilde\nabla_{\pi}h(F^h(A\pi))-\tilde\nabla_{\pi}h(\pi)\r)'\l(F^h(A\pi)-\pi\r)\,.$$
It follows from (\ref{dgammadt0}), (\ref{dgammadt1}), and (\ref{dgammadt2}), that, for $\eta < \eta^*$ and $t \geq t_3$,  
\be\label{dgammadt3}
\frac{\de}{\de t}\Gamma(t)\le-\eta\Upsilon(\pi(t))+M\eta^2\,.\ee 
From the strict convexity of $h(\pi)$ on the simplex $\Pi$, one finds that $\Upsilon(\pi)\ge0$ for all $\pi$, with equality iff $\pi=\pi^h$. Now, let  
$$\delta(x):=\l\{\ba{lcl}\sup\{||A\pi-f^h||:\,\Upsilon(\pi)\le Mx\}+M x &\se&0\le x<\eta^*\\
\tilde C \sqrt{|\mc E|}&\se& x\ge\eta^*\, ,\ea\r.$$
where $\tilde C:=\max\{\tilde C_e: e \in \mc E\}$, with $\tilde C_e$ as defined in Lemma~\ref{lemmaawayC2}.
It can be verified that $\delta(x)$ is right-continuous, nondecreasing, and such that $\delta(0)=0$. 
Then, (\ref{eq:f-Api}) and (\ref{dgammadt3}) imply that,  for $\eta < \eta^*$,
$$\limsup_{t\to+\infty}||f(t)-f^h||\le\delta(\eta)\,.$$
For $\eta\ge \eta^*$, the above is clearly true since $f(t) \in [0,\tilde C]^{\mc E}$ by Lemma~\ref{lemmaawayC2} and $f^h\in A\Pi\subseteq[0,1]^{\mc E}$.
Together with Lemma \ref{lemmahto0}, this completes the proof of Theorem \ref{thm:main-unperturbed}.

\section{Simulations}
\label{sec:sim}
In this section, we present results from numerical experiments. We performed several experiments with different graph topologies and for values of $\eta$ ranging from 0.01 to 100. In all the cases, we found that the trajectories converge exactly to the perturbed Wardrop equilibrium, i.e., $\delta(\eta)$ in Theorem~\ref{thm:main-unperturbed} was estimated to be uniformly zero. We suspect that this might be because of the exponential convergence also of the slow scale dynamics. Additionally, we compared the convergence of the trajectories corresponding to local decision function from Example \ref{ex:i-logit} with trajectories corresponding to local decision function of the form \be\label{triviallocaldecfunct} G^v_e\left(f_{\mc E^+_v},\pi \right)=f_e^{\pi}/\sum_{j\in\mc E^+_v}f_j^{\pi}\,,\qquad 
\forall f_{\mc E^+_v}\in\mc F_v,\ \forall e \in \edgeset_v^+\,.\ee
The latter corresponds to the case when the drivers do not take into account the local observation on the currently observed flow, and always act in a way that is consistent with their aggregate path preference. 
We found that the trajectories corresponding to local decision function in \eqref{triviallocaldecfunct} converged faster than the trajectories corresponding to the local decision function in Example~\ref{ex:i-logit}.

\begin{figure}
\centering
\includegraphics[width=0.85\linewidth]{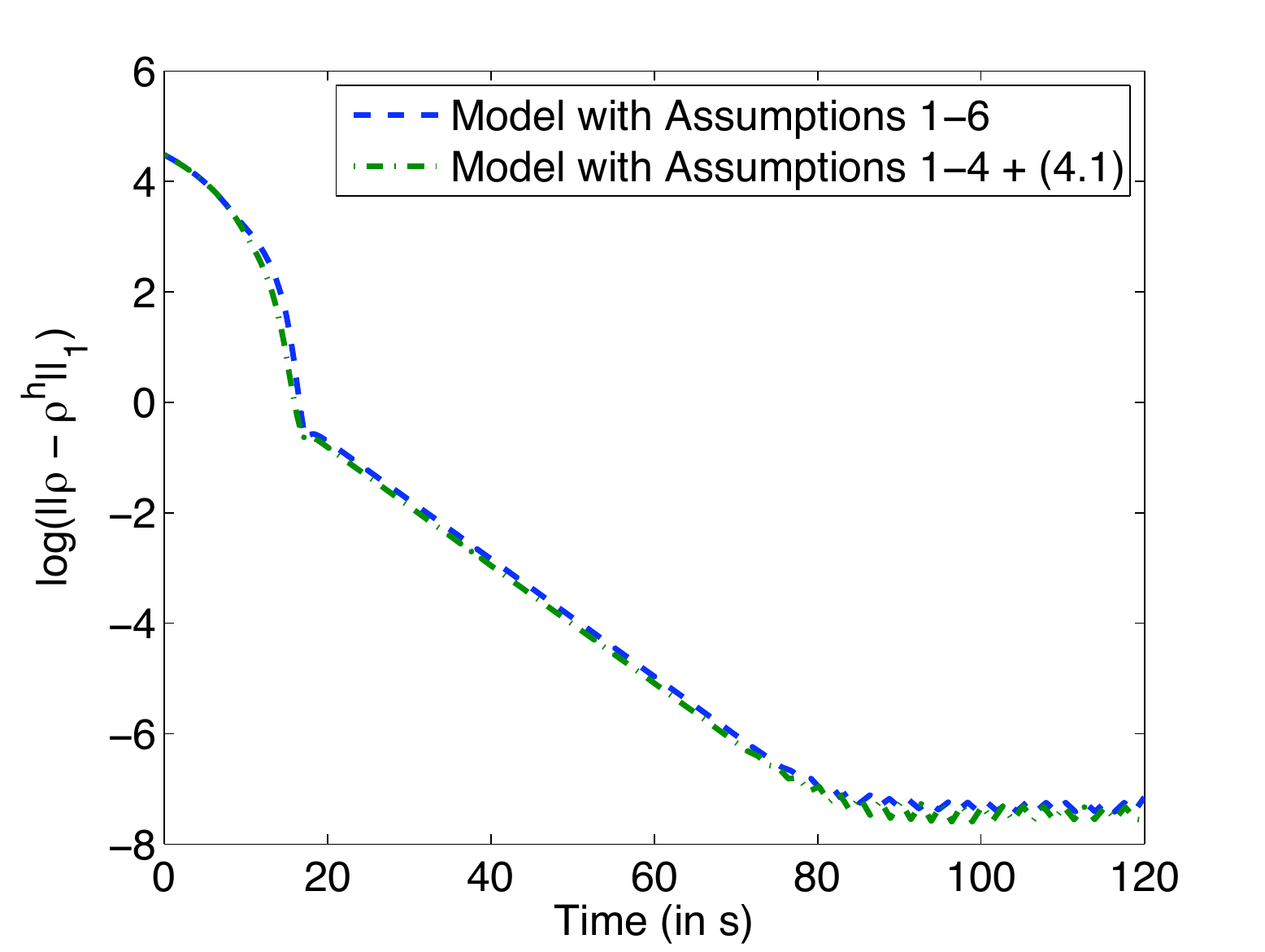}
\caption{Log-linear plot for comparison of the evolution of $\|\rho(t)-\rho^h\|_1$ for the local decision function of Example \ref{ex:i-logit} versus the local decision function of (\ref{triviallocaldecfunct}).}
\label{fig:comparisonplot}
\end{figure}

We demonstrate these findings through an illustrative example. For this example, the parameters were selected as follows:
\begin{itemize}
\item graph topology $\graph$ as shown in Figure~\ref{fig:graph},
\begin{figure}
\begin{center}
\scalebox{1.0}
{
\begin{tikzpicture}

\path (0,0) node(a) [circle,draw,fill=blue!50!white] {$0$}
		(2,1.5) node (b) [circle,draw,fill=blue!50!white] {$1$}
	(2, 0) node (c) [circle,draw,fill=blue!50!white] {$2$}
   (2, -1.5) node (d) [circle,draw,fill=blue!50!white] {$3$}
   (4, 1) node (e) [circle,draw,fill=blue!50!white] {$4$}
   (4, -1) node (f) [circle,draw,fill=blue!50!white] {$5$}
    (6, 2.5) node (g) [circle,draw,fill=blue!50!white] {$6$}
    (6, -2.5) node (h) [circle,draw,fill=blue!50!white] {$7$}
    (8, 0) node (i) [circle,draw,fill=blue!50!white] {$8$};
	
\draw[thick,->] (-1,0) -- (a); 
\draw[very thick,->] (a) -- (b);
\draw[very thick,->] (a) -- (c);
\draw[very thick,->] (a) -- (d);
\draw[very thick,->] (b) -- (e);
\draw[very thick,->] (c) -- (e);
\draw[very thick,->] (c) -- (f);
\draw[very thick,->] (d) -- (f);
\draw[very thick,->] (e) -- (i);
\draw[very thick,->] (f) -- (i);
\draw[very thick,->] (e) -- (g);
\draw[very thick,->] (f) -- (h);
\draw[very thick,->] (b) -- (g);
\draw[very thick,->] (d) -- (h);
\draw[very thick,->] (g) -- (i);
\draw[very thick,->] (h) -- (i);

\draw (-0.65,0.3) node {$1$};	
\draw (1,1.1) node {$e_1$};
\draw (1,0.2) node {$e_2$};	
\draw (1,-1.1) node {$e_3$};	
\draw (3.2,1.4) node {$e_4$};
\draw (2.8,0.7) node {$e_5$};	
\draw (3,-0.3) node {$e_6$};	
\draw (2.9,-1.1) node {$e_7$};	
\draw (3.8,-2.2) node {$e_8$};	
\draw (3.8,2.2) node {$e_9$};
\draw (5.2,1.6) node {$e_{10}$};	
\draw (5.2,-1.6) node {$e_{11}$};	
\draw (6,0.7) node {$e_{12}$};	
\draw (6,-0.7) node {$e_{13}$};	
\draw (7.2,1.5) node {$e_{14}$};	
\draw (7.2,-1.5) node {$e_{15}$};		
	
\end{tikzpicture}
}
\end{center}
\caption{The graph topology used in simulations.}
\label{fig:graph}
\end{figure}
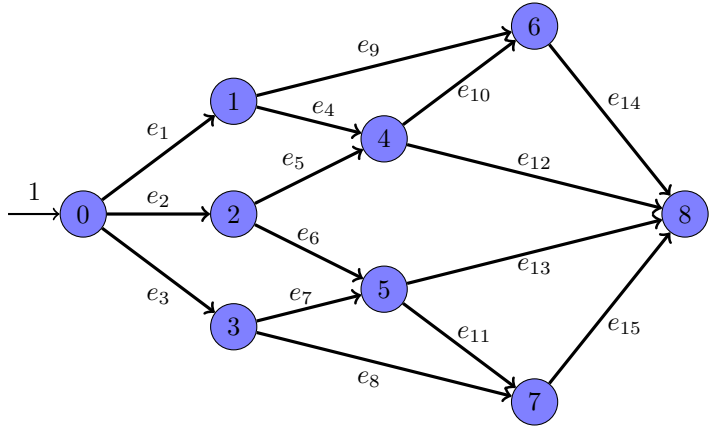
\item link-wise flow functions as given by \eqref{eq:flowfunc-example} with $C_1= 2$ and $\theta_e=1$, for all $e\in\mc E$;
\item $F^h$ as in \eqref{eq:logit} with $\beta=1$,
\item $G$ as in \eqref{eq:i-logit} with $\gamma=1$,
\item initial conditions: $\pi_e(0)=1/15$ for all $e\in\mc E$, $\rho_{e_1}(0)=\rho_{e_{12}}(0)=5$, $\rho_{e_2}(0)=\rho_{e_{6}}(0)=\rho_{e_{8}}(0)=7$, $\rho_{e_3}(0)=\rho_{e_{7}}(0)=3$, $\rho_{e_4}(0)=6$, $\rho_{e_5}(0)=1$, $\rho_{e_9}(0)=9$, $\rho_{e_{10}}(0)=10$, $\rho_{e_{13}}(0)=12$, $\rho_{e_{14}}(0)=4$,  $\rho_{e_{15}}(0)=8$. 
\item $\eta=0.1$.
\end{itemize}
For these values, $\rho^h:=\mu^{-1}(f^h)$ was numerically calculated by implementing a gradient descent algorithm for the potential function as given in \eqref{eq:potential-def}.
The evolution of the 1-norm distance of $\rho$ from $\rho^h$ is plotted on a log-linear scale in Figure~\ref{fig:comparisonplot} for two cases: (i) local route choice decision function of Example \ref{ex:i-logit}, and (ii) local decision function given in (\ref{triviallocaldecfunct}). Figure~\ref{fig:comparisonplot} also shows that there is no significant difference between the convergence of trajectory corresponding to local decision function in \eqref{triviallocaldecfunct} and the trajectory corresponding to the local decision function of Example~\ref{ex:i-logit}. However, as we increase $\eta$, we observed that the trajectory corresponding 
to the local decision function in \eqref{triviallocaldecfunct} converge faster than the trajectory corresponding to the local decision function of Example~\ref{ex:i-logit}. 


\section{Conclusion}
\label{sec:conclusion}
In this paper, we analyzed the stability of Wardrop equilibria in dynamical transportation networks characterized by dual temporal and spatial scales of the drivers' route choice behavior. We showed that, if the frequency of updates of  path preferences is sufficiently small, then the state of the transportation network ultimately approaches a neighborhood of the Wardrop equilibrium. The technical approach relied on establishing relevant properties for the resultant two time-scale dynamics independently using tools from evolutionary game dynamics, and cooperative dynamical systems, and then using singular perturbation techniques to establish sufficient conditions for the stability of the Wardrop equilibrium for the coupled system. Our results contribute to providing a stronger evidence in support of the significance of Wardrop's postulate of equilibrium for a transportation network. They may be read as a sort of robustness of such equilibrium notion with respect to non-persistent perturbations of the network.

There are several possible directions for future work. 
We plan to formally justify our dynamical model as a macroscopic approximation of the underlying driver level microscopic process.  We also plan to extend our analysis to the case with multiple origin-destination pairs and possibly cyclic topologies. We also plan to study the effect of persistent, and possibly adversarial, perturbations on the traffic dynamics under driver behavior model similar to the one considered in this paper, e.g., see ~\cite{Como.Savla.ea:MTNS10}.

{\small
  \bibliographystyle{siam}%
  \bibliography{KS-transportation}
}

\end{document}